\documentclass[11pt,a4paper]{article}
\usepackage[english]{babel}
\usepackage[cp1250]{inputenc}
\usepackage{hyperref}
\usepackage[OT1]{fontenc}
\usepackage{latexsym}
\usepackage{graphicx}

\makeatletter

\newtheorem{df}{Definition}[section]

\newcommand{\F}[0]{\mathcal{F}}
\newcommand{\T}[0]{\mathcal{T}}
\newcommand{\B}[0]{\mathcal{B}}

\author{{\large Anna Dudek$^{1,}$}\footnote{e-mail: aedudek@uci.agh.edu.pl}\quad {\large Maciej Goćwin$^{1,}$\footnote{e-mail: gocwin@uci.agh.edu.pl}}
 \quad {\large Jacek Leśkow$^2$}\\ {\footnotesize
 $^1$\it Department of Applied Mathematics, AGH University of Science
and Technology,}\\ {\footnotesize \it Al. Mickiewicza 30, 30-059
Cracow, Poland, phone: +48 12 617 4405}\\{\footnotesize $^2$\it
Econometrics Department, Nowy S\c{a}cz Graduate School of Business
- National Louis University,}\\ {\footnotesize \it Zielona 27,
33-300 Nowy S\c{a}cz, Poland, phone: +48 507 096 505}}

\title{Simultaneous confidence bands for the integrated hazard function}

\begin{document}

\maketitle
\begin{abstract} The construction of the simultaneous
confidence bands for the integrated hazard function is considered.
The Nelson--Aalen estimator is used. The simultaneous confidence
bands based on bootstrap methods are presented. Two methods of
construction of such confidence bands are proposed. The weird
bootstrap method is used for resampling. Simulations are made to
compare the actual coverage probability of the bootstrap and the
asymptotic simultaneous confidence bands. It is shown that the
equal--tailed bootstrap confidence band  has the coverage
probability  closest to the nominal one. We also present
application of our confidence bands to the data regarding
 survival after heart transplant.
\end{abstract}

\bigskip
\noindent {\bf Keywords:} bootstrap, intensity function,
multiplicative intensity model, Nelson--Aalen estimator, point
process, simultaneous confidence bands.

\bigskip
\noindent {\bf MSC:} 65C60, 62F25, 62F40.

\begin{section}{Introduction and summary}
In biomedical settings, the multiplicative intensity model
introduced by Aalen has many applications. This is a model for
point processes observed on a fixed time interval for which the
stochastic intensity is decomposed into deterministic function
$\alpha(t)$ and stochastic process $Y(t)$. The $\alpha(t)$
function may be considered as an individual force of transition at
time $t$ and $Y(t)$ as a number at risk just before time $t$.
\\ \\
In broad terms what makes survival data special is the presence of
censored data. To analyze such data by the multiplicative
intensity model a general assumption of independent censoring is
required, which means that at any time $t$ the survival experience
in the future is not statistically altered by censoring and
survival experience in the past. The censoring mechanism is
modelled by $Y$ process and has not any influence on the $\alpha$
function.
\\ \\
In the survival analysis the most interesting is to estimate the
survivor function and the integrated hazard function. In this
paper we consider the latter, which is estimated by the
Nelson--Aalen estimator. An interpretation of this  estimator is
difficult without construction of some confidence intervals. From
our perspective, the pointwise intervals are not totaly
satisfactory while one wants to construct  confidence region for
the whole curve simultaneously for all points.
\\ \\
The construction of the simultaneous confidence bands is difficult
since we need the uniform consistency property. However, such
confidence bands are badly needed in practical applications. For
example, in the works related with ours like in the papers of
Cowling, Hall , Phillips (\cite{hall}) and Snethlage
(\cite{sneth}) but also in the time series analysis (Leśkow and
Wronka \cite{wronka}) and the nonparametric regression (Loader
\cite{loader}).
\\ \\
The formula for the asymptotic confidence interval for the
Nelson--Aalen estimator is known, however, it is very complicated
and does not work well for small samples (see \cite{gill}). An
alternative approach is through the use of bootstrap methods. This
idea was first introduced by Efron (\cite{efron}) and later
developed in many papers (also in cited above). Bootstrapping of
the point processes is not yet fully explored. Some results are
presented in \cite{braun1998} and \cite{braun2003}. The Poisson
process context is treated in the paper \cite{hall}, however these
methods cannot be easily adapted to the multiplicative intensity
model.
\\ \\
The aim of our work is the construction of the bootstrap
simultaneous confidence bands for the Nelson--Aalen estimator. We
want to compare constructed bootstrap regions with the asymptotic
ones. We make simulations to check if the actual coverage
probability is close to a nominal one. In our calculations we use
the weird bootstrap method.
\\ \\
We show that for the small samples the bootstrap models have much
better coverage probabilities. Not only the actual coverage
probabilities of the bootstrap simultaneous confidence bands are
very close to nominal ones but also the left- and right--tail
error probabilities are almost equal.
\\ \\
Our paper is organized in the following way. Section
\ref{formulation} contains a short survey of basic results related
to the Nelson--Aalen estimator and the bootstrap for point
processes. Section \ref{simband} is dedicated to construction of
simultaneous confidence bands for the estimator considered. A
practical example related to heart transplant study is included in
Section \ref{example}, while Section \ref{results} contains
additional numerical results. Conclusions and open questions are
presented in Section \ref{conclusions}.
\end{section}

\begin{section}{Problem formulation}\label{formulation}
In our paper we construct the bootstrap simultaneous confidence
bands for the integrated intensity function. We use the weird
bootstrap introduced in \cite{gill}. We compare our results with
those presented in \cite{gill} and \cite{bie}. Application of
bootstrap is well motivated in the small sample case and when
censoring mechanism is quite complex. Moreover, the standard
asymptotic theory provides confidence intervals that are quite
difficult to apply. To construct bootstrap simultaneous confidence
bands we applied one of the methods proposed in \cite{hall}.
\\ \\
While defining our problem we follow \cite{gill} (page 176). We
consider a continuous--time interval $\T$ which may be of the form
$[0,\tau]$ or $[0,\tau)$ for a given terminal time $\tau$,
$0<\tau\leq\infty$. Let $(\Omega,\F)$ be a measurable space
equipped with a filtration $(\F_{t}, t\in\T)$. We define on
$(\Omega,\F)$ a counting process $\textbf{N}=(N(t),t\in\T)$
adapted to the filtration such that its stochastic intensity
function $\lambda$ is of the form $\lambda(t)=\alpha(t)Y(t)$,
where $\alpha$ is nonnegative deterministic function and $Y$ is a
predictable process. For example, we can consider an initial group
$Y_{0}$ of patients with cancer after some medical treatment.
Although the patients enter the study at different calendar times,
we observe only their time since operation. In this case $\alpha(t)$
is the individual intensity of death and $Y(t)$ is the number at
risk at the moment of time $t$ e.g. number of patients who lived
till time $t$. For a practical example see Section \ref{example}.
\\
The only assumption we have to make about $\alpha$ is its
integrability,
\[
\int_{0}^{t}\alpha(s)ds<\infty \ \ \ \ \textrm{for all } t\in\T.
\]
We consider the Nelson--Aalen estimator $\widehat{A}$ for
\[
A(t)=\int_{0}^{t}\alpha(s)ds
\]
which is of the form
\[
\widehat{A}(t)=\sum_{j:T_{j}\leq t}\frac{1}{Y(T_{j})},
\]
where $T_{j}$ are jump times.
\\
We define an estimator for the mean squared error function as
\[
\widehat{\sigma}^{2}(t)=\sum_{j:T_{j}\leq t}\frac{Y(T_{j})-\Delta
N(T_{j})}{Y^{3}(T_{j})},
\]
where $\Delta N(T_{j})=N(T_{j})-N(T_{j-1})$.
\\ \\
Under the suitable assumptions the Nelson--Aalen estimator is
uniformly consistent on compact intervals  (see \cite{gill} page
190), which means:
\[
\sup_{s\in[0,t]}|\widehat{A}^{(n)}(s)-A(s)|\stackrel{p}\longrightarrow
0 \quad \textrm{as} \ n\rightarrow\infty \quad \textrm{for} \
t\in\T.
\]
The asymptotic distribution of the Nelson--Aalen estimator can be
obtained from Rebolledo's martingale central limit theorem (for
details see \cite{gill} page 190). It should be pointed out that
the problem of constructing simultaneous confidence bands requires
a version of the functional central limit theorem for the
cumulative intensity function. Such results can be found in
\cite{gill} (page 263), however the limiting distribution is quite
difficult to apply in practice. Moreover, it is still unknown what
form of the functional central limit theorem can be established
for $\alpha$ alone. (See also Section \ref{conclusions} for
additional remarks regarding this problem).
\\ \\
The results above can be used to construct pointwise confidence
intervals and simultaneous confidence bands for $A(t)$
(\cite{gill}). Unfortunately, formulae for the asymptotic
distributions are very complicated. That is why we want to apply
bootstrap methods to construct simultaneous confidence bands.
Bootstrapping of counting processes is not easy because such
processes are not based on i.i.d. samples. The problem is complex
and, thus, the methods for the general case are not known. There
are some results for the Poisson processes (see \cite{hall}),
however in this case one may get similar results without
simulations (see \cite{sneth}). Some methods of bootstrapping
point processes are also presented in \cite{braun1998} and
\cite{braun2003}.
\\ \\
In our paper we apply the weird bootstrap method. The idea is
based on the fact that the asymptotic distribution of
$a_{n}(\widehat{A}-A)$ has independent increments and
$Var(d\widehat{A}(t)|\F_{t-})=dA(t)(1-dA(t))/Y(t)$. The following
definition is quoted from \cite{gill}.
\begin{df}{The Weird Bootstrap}
\\Given $N$, $Y$, and $\widehat{A}$, let $N^{*}$ be a process with
independent binomial $(Y(t),\Delta\widehat{A}(t))$ distributed
increments at the jump times of $N$, constant between jump times.
Let $\widehat{A}^{*}=\int dN^{*}/Y$. Estimate the distribution of
$\widehat{A}-A$ by the conditional distribution, given $N$ and
$Y$, of $\widehat{A}^{*}-\widehat{A}$.
\end{df}
For the proof of consistency of this method see \cite{gill} (page
220).
\\
The word {\it weird} is not accidental. In every time point
$t\in\T$ every individual at risk from the set $Y(t)$ has the same
probability of a failure. However, the event at the time $t$ does
not exert any influence on any other time moment $s\in\T$.
\\ \\
The problem of bootstrapping point processes is not completely
solved and quite challenging. Some partial solution are discussed
in \cite{hall}, \cite{hinkley}, \cite{braun1998} and
\cite{braun2003}. In the next section we use this method of
bootstrapping to construct the simultaneous confidence bands.

\end{section}

\begin{section}{Simultaneous confidence bands}\label{simband}
The Nelson--Aalen point estimator is difficult to interpret
without some idea of its accuracy. Resolving this problem requires
constructing confidence intervals or confidence bands. These bands are also
quite interesting because of their hypothesis testing
interpretation. We can think of confidence bands as a one--sample
test statistics with a null hypothesis $A=A_{0}$ which is rejected
at significance level $\theta$ if $A_{0}$ is not completely
contained in the band. In this case pointwise confidence intervals
are not satisfactory. That is why we introduce simultaneous
confidence bands.
\begin{df}{Confidence region}
\\Let $\B$ denote a connected, nonempty, random subset of the
rectangle $[0,\tau]\times[0,\infty)$, such that $\B\cap\{(x,y):
0\leq y< \infty\}$ is nonempty for each $x\in[0,\tau]$. We call
$\B$ a confidence region for $A$ over the set $S\in[0,\tau]$ with
a coverage probability $(1-\theta)$ if $P\{(x,A(x))\in\B \ for\
all\ x\in S )\}=\theta$.
\end{df}
In our paper $S$ is always an interval.
\\ \\
Simultaneous confidence bands may be constructed in many different
ways. The authors of the book \cite{gill} (page 209) proposed two
types of such bands: EP--band (equal precision band) and HW--band
(Hall--Wellner band). These confidence bands are based on the
asymptotic distribution of the Nelson--Aalen estimator on compact
intervals which can be derived  from the martingale central limit theorem.
\\
Both EP- and HW--band for $A$ on $[t_{1},t_{2}]$ are of the form
\[
\widehat{A}(s)\pm
a_{n}^{-1}K_{q,\theta}(c_{1},c_{2})(1+a_{n}^{2}\widehat{\sigma}^{2}(s))/q(\frac{a_{n}^{2}\widehat{\sigma}^{2}(s)}{1+a_{n}^{2}\widehat{\sigma}^{2}(s)})
\]
with $K_{q,\theta}(c_{1},c_{2})$ being the upper percentile of the
distribution of
\[
\sup_{x\in[c_{1},c_{2}]}|q(x)W^{0}(x)|,
\]
where $W^{0}$ denotes the standard Brownian bridge.
\\
The constants $c_{1}$ and $c_{2}$ can be approximated by
\[
\widehat{c}_{i}=\frac{a_{n}^{2}\widehat{\sigma}^{2}(t_{i})}{1+a_{n}^{2}\widehat{\sigma}^{2}(t_{i})},
\]
where $a_{n}=\sqrt n$ is a normalizing factor and $n$ is the
number of individuals at study. \\
For EP--band $q$ is chosen as $q_{1}(x)=\{x(1-x)\}^{-1/2}$ which
yields the confidence bands proportional to the pointwise ones.
For HW--band $q$ is chosen as $q_{2}(x)=1$.
\\
In both cases $\theta$ percentile of the asymptotic distribution
are difficult to obtain. These bands also perform badly even with
the sample size of 100--200 \cite{bie}. Because of this reason one
may consider some transformations to improve the approximation to
the asymptotic distribution \cite{gill} (page 211).
\\ \\
To avoid such problems we consider bootstrap simultaneous
confidence bands. The authors of the paper \cite{hall} proposed a
few different methods of constructing these bands. In our
calculations we use the weird bootstrap method. Our construction
of bootstrap--t confidence regions for $A$ is based on the
bootstrap approximation
\[
T^{*}(x)=\frac{\widehat{A}^{*}(x)-\widehat{A}(x)}{\widehat{\sigma}(x)},\
\ x\in\T
\]
of
\[
T(x)=\frac{\widehat{A}(x)-A(x)}{\widehat{\sigma}(x)},\ \ x\in\T.
\]
For details see \cite{gill}.
\\ \\
Below we present two bootstrap confidence bands:

\def\theenumi{ \bf \arabic{enumi}}

\begin{enumerate}
\item Confidence region is defined by
\[
\B_{1}=\{(x,y):x\in S,\
\max[0,\widehat{A}(x)-t_{1}\widehat{\sigma}(x)]\leq
y\leq\widehat{A}(x)+t_{1}\widehat{\sigma}(x)\},
\]
where $t_{1}$ is chosen such that
\[
P\{|T^{*}(x)|\leq t_{1},\ \textrm{all}\ x\in S|N,Y\}=1-\theta.
\]
The main feature of this region is that at the point $x$ its width
is proportional to $\widehat{\sigma}(x)$.

\item In many applications populations cannot be modelled via
symmetric distributions. The only reasonable choice is a strongly
skewed distribution. In all of the previous presented intervals,
skewness was not taken into consideration. This has a quite
negative impact on the coverage probability. To adjust for
skewness of the distribution one could construct a region which
the left- and right--tail error probabilities are equal. This kind
of the region is of the form
\[
\B_{2}=\{(x,y):x\in S,\
\max[0,\widehat{A}(x)-t_{3}\widehat{\sigma}(x)]\leq
y\leq\widehat{A}(x)-t_{2}\widehat{\sigma}(x)\},
\]
where $t_{2}$ and $t_{3}$ are chosen such that
\[
P\{t_{2}\leq T^{*}(x)\leq t_{3},\ \textrm{all}\ x\in
S|N,Y\}=1-\theta
\]
and
\[
P\{T^{*}(x)\leq t_{2},\ \textrm{all}\ x\in S|N,Y\}=P\{T^{*}(x)\geq
t_{3},\ \textrm{all}\ x\in S|N,Y\}.
\]
\end{enumerate}
In the next section we present an example of applying such bands.

\end{section}

\begin{section}{Practical example}\label{example}
We take under the consideration the group of 64 patients after
heart transplant. The data we use are taken from
\cite{kalbfleisch}, Appendix A, pages 387-389. In our approach,
the risk is defined as the rejection of the transplant so the time
between the operation and the rejection is considered. 35
observations are censored. The censoring was present if patients
were alive at the end of the study or lost to follow--up. The 95\%
confidence bands simultaneous with respect to the time argument
were constructed in the time bandwidth between day 20 and day 1200
of the observation. The construction of such confidence interval
was based on Nelson--Aalen estimator.
\begin{figure}[p]
\begin{center}
\includegraphics[width=\textwidth]{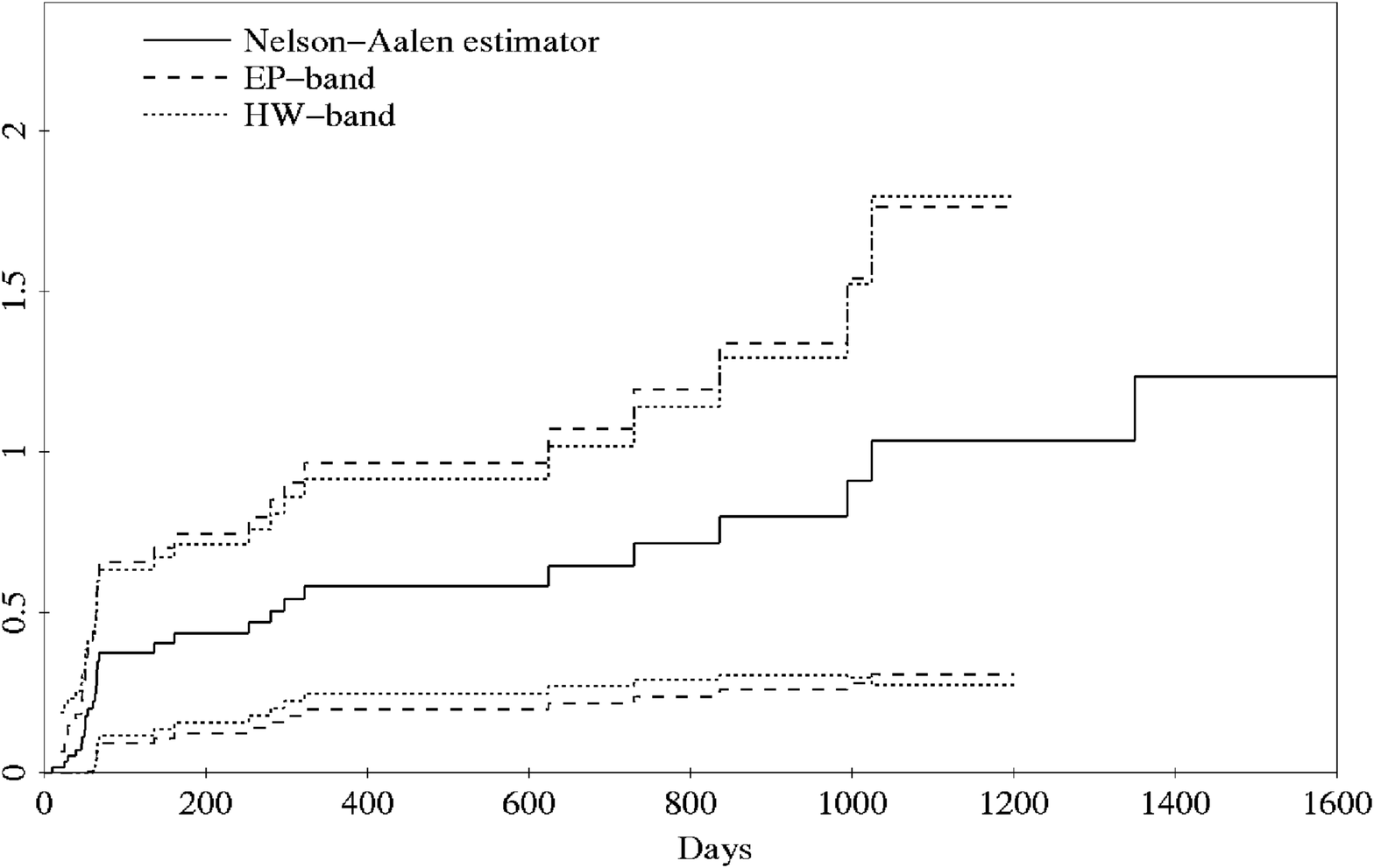}
\caption{HW- and EP--band}\label{fig1}
\includegraphics[width=\textwidth]{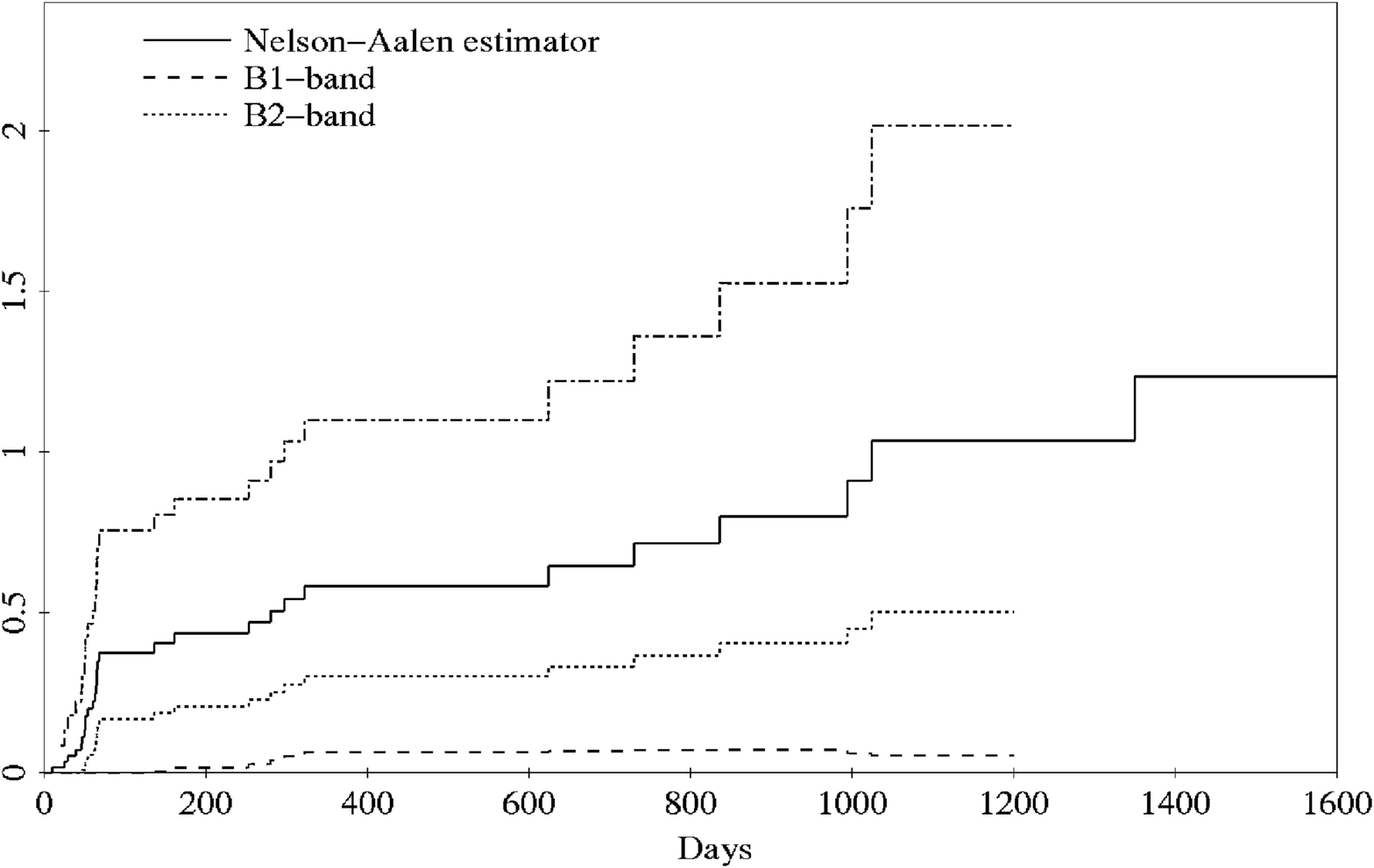}
\caption{$\B1$- and $\B2$--band}\label{fig2}
\end{center}
\end{figure}
Figure \ref{fig1} presents the Nelson--Aalen estimator together
with HW and EP bands and Figure \ref{fig2} with $\B1$ and $\B2$
bootstrap simultaneous confidence bands. Note that $\B1$, EP and
HW bands are symmetric. Only $\B2$ is not symmetric. The upper
bands of $\B1$ and $\B2$ are covering themselves. The lower band
of $\B1$ is noticeably too low. It suggests that $\B1$ is too
wide. HW and EP bands are close to each other but EP is
significantly broader during the most part of the time interval.
Moreover, $\B2$ is shifted upwards compared
with the asymptotic simultaneous confidence bands. \\
Now we will verify the actual coverage probability for the
considered bands.
\end{section}

\begin{section}{Numerical results}\label{results}
Our aim is to compare the coverage probability for asymptotic and
bootstrap simultaneous confidence bands. Our simulations are based
on the multiplicative model for the intensity function
$\lambda(t)=Y(t)\alpha(t)$. We concentrate on a few typical
examples of the $\alpha$ function. To generate process Y we first
choose the beginning value $Y_{0}$ (the number of individuals at
risk) and next for every individual the time of termination is
sampled from exponential distribution with the mean value $0.25$.
Having such $Y$ we generate the underlying point process.
\\
For our study we chose four functions:
\begin{eqnarray*}
 \alpha_{1}(t)&=&\frac{5}{3}, \\
 \alpha_{2}(t)&=& \frac{5}{6}+10(t-0.5)^{2},\\
 \alpha_{3}(t)&=&\frac{5}{3}+10(t-0.5)^{3}, \\
 \alpha_{4}(t)&=&2.5-10(t-0.5)^{2}.
\end{eqnarray*}
Curves of such kinds can be applied in biomedicine, insurance and
demography. For example the U--shaped functions may reflect
behavior of the intensity of death and the inverted U--shaped
functions may describe the intensity of birth. These shapes are
reflected in the equation of $\alpha_{2}$ and $\alpha_{4}$
functions. Figure \ref{funkcje} shows these intensity functions
and Figure \ref{calki} presents integrated versions of these
functions.
\begin{figure}[p]
\begin{center}
\includegraphics[width=\textwidth]{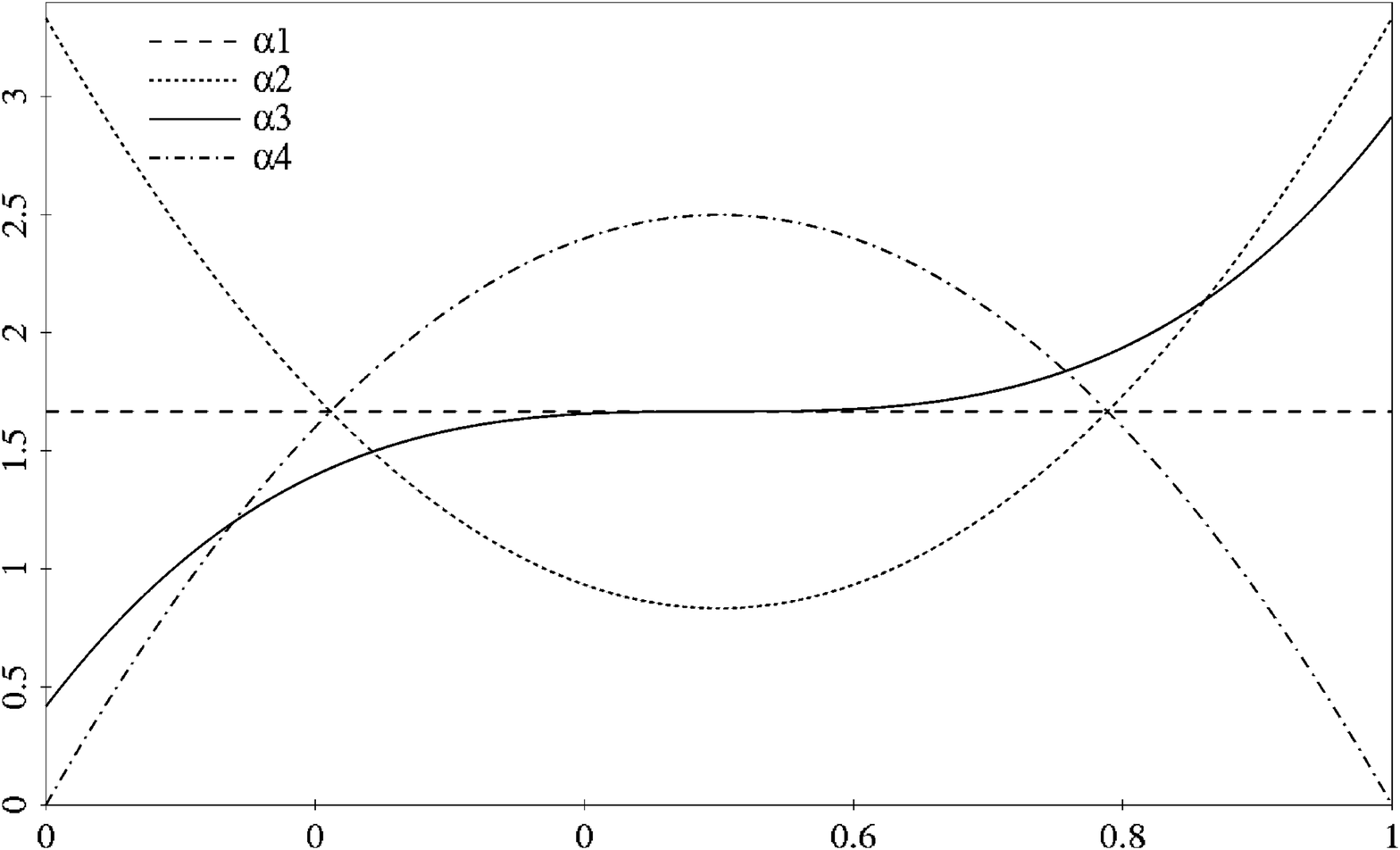}
\caption{Intensity functions}\label{funkcje}
\includegraphics[width=\textwidth]{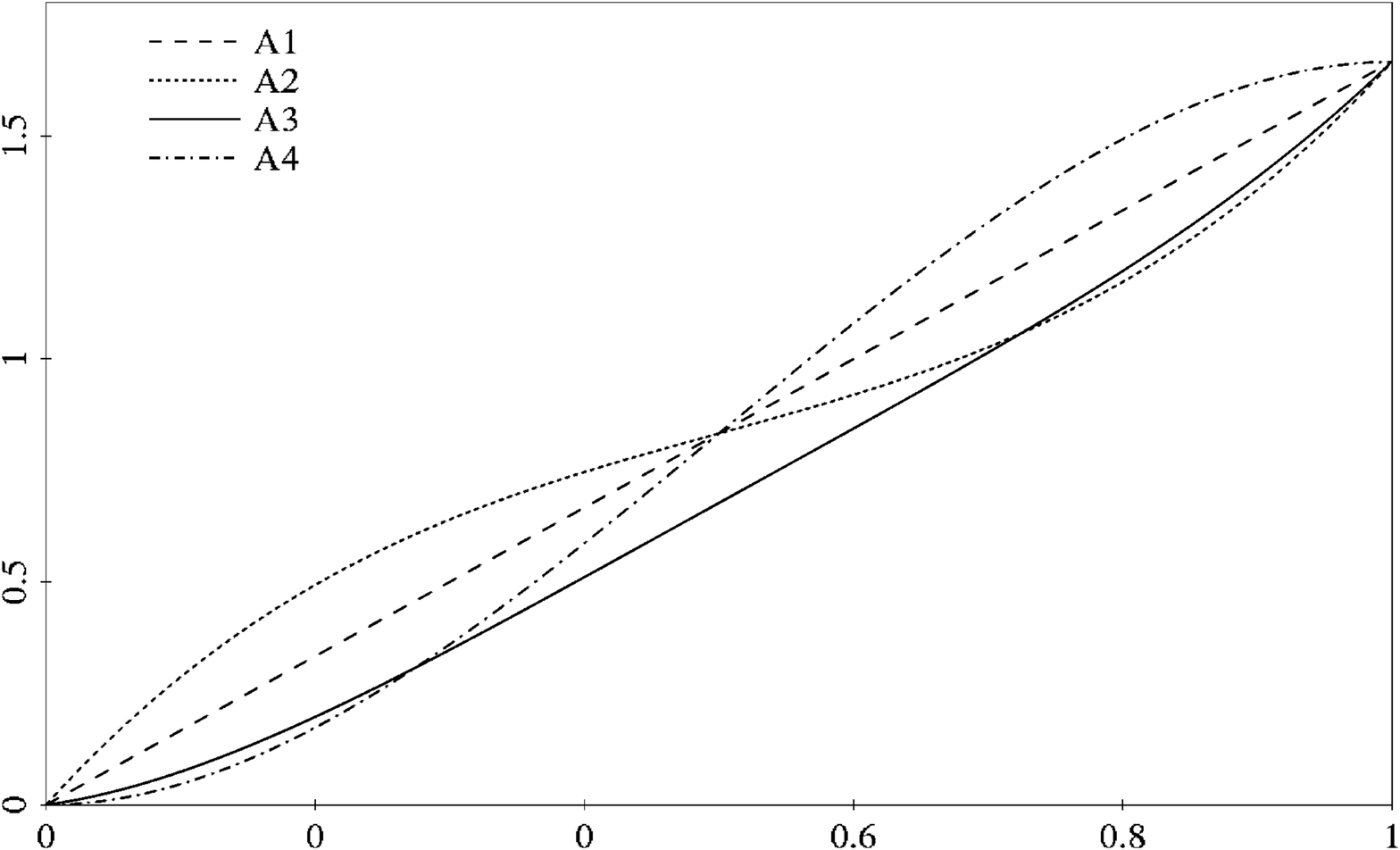}
\caption{Integrated intensity functions}\label{calki}
\end{center}
\end{figure}

 \ \\
 We make simulations for the interval $S=[0.2;0.8]$, the number of
 bootstrap resamples $B=200$ and initial number at risk
 $Y_{0}:25,50,75$. In Table \ref{tab1} we show the actual coverage
 probability, when the nominal coverage
 probability is 0.95 and the number of iterations is equal to 10000.
 For every $Y_{0}$, $\alpha_{i}$ function ($i=1 \dots 4$)
 and method of construction of the confidence
 region, the first and the second number in Table \ref{tab1} are the
 left- and right--tail error probabilities and the third is the
 actual coverage probability (all probabilities are measured in
 percentage).
\begin{table}[h]
\begin{center}
\begin{tabular}{|c|c|ccc|ccc|ccc|}\hline
 & & \multicolumn{9}{|c|}{$Y_{0}$} \\ \hline
 function & method & \multicolumn{3}{|c|}{25} & \multicolumn{3}{|c|}{50} & \multicolumn{3}{|c|}{75} \\ \hline
              & HW & 0.4 & 11.5 & 88.1 & 0.6 & 5.7 & 93.7 & 0.8 & 4.1 & 95.1 \\
  $\alpha_{1}$ & EP & 0.4 & 12.7 & 87.0 & 0.6 & 6.6 & 92.8 & 0.8 & 4.7 & 94.5 \\
              & $\B1$ & 0.0 & 4.0 & 96.0 & 0.1 & 2.5 & 97.4 & 0.3 & 2.2 & 97.5 \\
              & $\B2$ & 2.4 & 4.0 & 93.6 & 3.2 & 2.6 & 94.3 & 3.1 & 2.2 & 94.7 \\ \hline
              & HW & 0.5 & 10.1 & 89.5 & 0.9 & 6.5 & 92.7 & 1.0 & 4.7 & 94.3 \\
  $\alpha_{2}$ & EP & 0.5 & 10.0 & 89.5 & 0.8 & 6.6 & 92.6 & 1.0 & 5.0 & 94.1\\
              & $\B1$ & 0.1 & 3.6 & 96.3 & 0.2 & 3.4 & 96.4 & 0.3 & 2.6 & 97.1\\
              & $\B2$ & 2.6 & 3.6 & 93.8 & 3.0 & 3.6 & 93.4 & 3.1 & 2.9 & 94.1\\ \hline
              & HW & 0.2 & 13.2 & 86.6 & 0.7 & 9.8 & 89.5 & 0.8 & 4.7 & 94.5 \\
  $\alpha_{3}$ & EP & 0.2 & 14.8 & 85.0 & 0.5 & 9.4 & 90.1 & 0.7 & 5.6 & 93.7 \\
              & $\B1$ & 0.0 & 4.7 & 95.3 & 0.9 & 2.9 & 96.3 & 0.1 & 1.6 & 98.3 \\
              & $\B2$ & 2.3 & 4.7 & 93.0 & 2.9 & 2.2 & 95.0 & 2.9 & 1.7 & 95.5 \\ \hline
              & HW & 0.7 & 13.8 & 85.8 & 0.5 & 6.7 & 92.8 & 0.9 & 4.4 & 94.8 \\
 $\alpha_{4}$ & EP & 1.0 & 16.6 & 83.1 & 0.3 & 8.5 & 91.2 & 0.8 & 5.3 & 93.9 \\
              & $\B1$ & 1.0 & 5.2 & 94.7 & 0.1 & 2.3 & 97.6 & 0.1 & 1.5 & 98.4 \\
              & $\B2$ & 2.6 & 5.1 & 92.5 & 2.7 & 2.3 & 95.0 & 3.3 & 1.5 & 95.2 \\
\hline
\end{tabular}
\end{center}
\caption{Actual coverage probability}\label{tab1}
\end{table}

As we expected, HW- and EP--band perform quite badly for the small
samples. Especially for $Y_{0}=25$ the actual coverage probability
is $5\%$ to $10\%$ less then it should be. This happens because
these are asymptotic bands and in our case the number of jumps of
the point processes is not big enough to apply the asymptotic
distribution. For $Y_{0}=50$ the actual coverage probability for
these bands is better but always remains about $3\%$ smaller than
the nominal one. For $Y_{0}=75$ all results are satisfactory. The
first of the bootstrap confidence intervals which we proposed
performs well for small $Y_{0}$ but when the number of jumps rises
it remains consistently too wide. The equal--tailed bootstrap
confidence band ($\B2$) behaves well in all considered situations.
Its actual coverage probability is always close to nominal, even
in the case of small beginning number at risk (when the asymptotic
bands fail). Our simulations also show that the left--side failure
probability for the EP- and HW--band is significantly too small.
Its value is below $1\%$. This means that our functions
$\alpha_{i}(t)$ almost never cross the lower band of the
confidence region e.g. the lower band goes too far away from the
estimator. The advantage of the $\B2$ region is the equal tailed
feature. The lack of coverage probabilities for the left--hand
case and the right--hand case are almost equal.
\\
We checked empirically that $\B2$ is the optimal choice.
Independently of the beginning number at risk it has a coverage
probability close to the nominal one and, what is very important,
it insures almost equally divided failure probability.
\\ \\
Now we compare our results with those presented in \cite{bie}. The
authors of \cite{bie} proposed arcsine- and logarithmic--transform
of the Nelson--Aalen estimator. They considered the modifications
of EP- and HW--band which use these transformations. Such
constructed asymptotic simultaneous confidence bands perform
satisfactionary for sample size as low as 25. \\
Using simulation methods presented before we compare the behaviour
of these bands to the bootstrap band $\B2$. The results are
presented in Table \ref{tab2}. AHW and AEP denote the
arcsine--trasform of HW- and EP--band respectively. The
logarithmic--transform bands are denoted by LHW and LEP.
\\ \\
\begin{table}[h]
\begin{center}
\begin{tabular}{|c|c|ccc|ccc|ccc|}\hline
 & & \multicolumn{9}{|c|}{$Y_{0}$} \\ \hline
 function & method & \multicolumn{3}{|c|}{25} & \multicolumn{3}{|c|}{50} & \multicolumn{3}{|c|}{75} \\ \hline
              & AHW & 2.2 & 6.6 & 91.2 & 2.7 & 3.1 & 94.2 & 2.7 & 2.0 & 95.3 \\
              & AEP & 2.3 & 5.9 & 91.8 & 2.6 & 3.3 & 94.1 & 2.7 & 2.0 & 95.3 \\
 $\alpha_{1}$ & LHW & 3.1 & 5.1 & 91.8 & 3.5 & 2.4 & 94.1 & 3.4 & 1.4 & 95.2 \\
              & LEP & 3.2 & 3.8 & 93.0 & 3.4 & 2.0 & 94.6 & 3.3 & 1.3 & 95.4 \\
              & $\B2$ & 2.4 & 4.0 & 93.6 & 3.2 & 2.6 & 94.3 & 3.1 & 2.2 & 94.7 \\ \hline
              & AHW & 2.3 & 4.8 & 92.9 & 2.4 & 3.2 & 94.4 & 2.4 & 2.2 & 95.4 \\
              & AEP & 2.3 & 4.5 & 93.2 & 2.3 & 3.3 & 94.4 & 2.4 & 2.3 & 95.3\\
 $\alpha_{2}$ & LHW & 3.0 & 3.2 & 93.8 & 3.0 & 2.4 & 94.6 & 2.9 & 1.6 & 95.5\\
              & LEP & 3.1 & 2.4 & 94.5 & 2.9 & 2.0 & 95.1 & 3.0 & 1.5 & 95.5 \\
              & $\B2$ & 2.6 & 3.6 & 93.8 & 3.0 & 3.6 & 93.4 & 3.1 & 2.9 & 94.1\\ \hline
              & AHW & 2.3 & 8.7 & 89.0 & 2.5 & 3.8 & 93.7 & 2.4 & 2.3 & 95.3 \\
              & AEP & 2.2 & 6.6 & 91.2 & 2.2 & 3.1 & 94.7 & 2.3 & 2.3 & 95.4 \\
 $\alpha_{3}$ & LHW & 3.8 & 6.2 & 90.0 & 3.5 & 2.8 & 93.7 & 3.6 & 1.8 & 94.6\\
              & LEP & 3.8 & 4.5 & 91.7 & 3.5 & 1.8 & 94.7 & 3.6 & 1.1 & 95.3 \\
              & $\B2$ & 2.3 & 4.7 & 93.0 & 2.9 & 2.2 & 95.0 & 2.9 & 1.7 & 95.5 \\ \hline
              & AHW & 2.3 & 8.8 & 88.9 & 2.5 & 3.0 & 94.5 & 2.5 & 2.3 & 95.2 \\
              & AEP & 2.3 & 7.3 & 90.4 & 2.3 & 3.1 & 94.6 & 2.4 & 2.1 & 95.5 \\
 $\alpha_{4}$ & LHW & 4.0 & 7.4 & 88.6 & 3.9 & 2.5 & 93.4 & 3.5 & 1.6 & 94.9\\
              & LEP & 3.6 & 5.7 & 90.7 & 3.6 & 1.7 & 94.7 & 3.3 & 1.1 & 95.6 \\
              & $\B2$ & 2.6 & 5.1 & 92.5 & 2.7 & 2.3 & 95.0 & 3.3 & 1.5 & 95.2 \\
\hline
\end{tabular}
\end{center}
\caption{Actual coverage probability}\label{tab2}
\end{table}
As might be expected for the sample size 50 and 75 all methods
give satisfactory results. For a sample size 25 the bootstrap
simultaneous confidence band $\B2$ has better coverage properties
than transformed asymptotic ones. The actual coverage probability
of $\B2$ is about $92.5\%$ for all $\alpha_{i}$ functions. It is
about $2\%$ closer to the nominal than the actual coverage
probability of the transformed bands. At first sight LEP seems to
be good choice but as the sample size grows it gets too wide.\\
However, considered transformations improve the actual coverage
probability and the left- and right--tail error probabilities of
the asymptotic bands $\B2$ is still the best choice.

\end{section}

\begin{section}{Conclusions}\label{conclusions}
In many applications , the hazard function is much more
interesting and relevant to estimate than the integrated hazard
function, but it is also more challenging to estimate. There are
several approaches to that problem, the histogram based sieve
estimator considered in Leśkow, Różański \cite{leskow} and Leśkow
\cite{leskow1988} being one of them. Unfortunately, the version of
functional central limit theorem of such estimator is still an
open question. Without such result construction of the simultaneous
confidence bands is impossible.
\\ \\
In our paper we showed that for the small samples the bootstrap
simultaneous confidence bands behave better than the asymptotic
ones. They also have better actual coverage probability. An
advantage of the equal--tailed type confidence region is the
balance of the left- and right--tail error probability. A
disadvantage of all simultaneous regions considered in this paper
is the lack of taking into a consideration the shape of the estimated
function. The integrated hazard function is always nondecreasing.
Unfortunately, the lower confidence band decreases sometimes. It
may be interesting to construct regions taking into consideration
the known features of the estimated function (for example monotonicity, unimodality). \\
The other curious problem is bootstrapping of the point process.
We consider only one method (the weird bootstrap). In the paper
\cite{hall} other methods are proposed but only for Poisson
processes. A method for obtaining bootstrap replicates for the
one--dimensional point process is presented in \cite{braun1998}
and its multi--dimensional version is also proposed. Because of
deficient coverage properties in some cases, Braun and Kulperger
proposed in \cite{braun2003} a technique for one--dimensional
point process which uses the idea of re--colouring presented in
\cite{hinkley}. It remains an open question if these methods can
be applied in a general case.
\end{section}


\begin{thebibliography}{00}
\bibitem{gill}
Andersen P.K., Borgan O., Gill R.D., Keiding N. (1992) {\it
Statistical Models Based on Counting Procces}. Springer.
\bibitem{bie}
Bie O., Borgan \O., Liest\o l K. (1987) {Confidence intervals and
confidence bands for the cumulative hazard rate function and their
small sample properties}. {\it Scand. J. Statist.} 14: 221-233.
\bibitem{Billingsley1968}
Billingsley P. (1968) {\it Convergence of Probability Measures}.
New York: Wiley.
\bibitem{braun1998}
Braun W.J., Kulperger R.J. (1998) {A bootstrap for point
processes}. {\it J. Statist. Comput. Simul.} 60: 129-155.
\bibitem{braun2003}
Braun W.J., Kulperger R.J. (2003) {Re-colouring the
Intensity-Based Bootstrap for point Processes}. {\it
Communications in Statistics Simulation and Computation} 32(2):
475-488.
\bibitem{hall}
Cowling A., Hall P., Phillips M.J. (1996) { Bootstrap confidence
regions for the intensity of a Poisson point process}. {\it Am.
Stat. Ass.} 91(436): 1516-1524.
\bibitem{hinkley}
Davison A.C., Hinkley D.V. (1999) {\it Bootstrap Methods and their
Applications}. Cambridge University Press.
\bibitem{efron}
Efron B., Tibshirani R.J. (1993) {\it An Introduction to the
Bootstrap}. Chapman\&Hall/CRC.
\bibitem{kalbfleisch}
Kalbfleisch J.D., Prentice R.L. (2002) {\it The Statistical
Analysis of Failure Time Data}. John Willey \& Sons, Inc.
\bibitem{leskow1988}
Leśkow J. (1988) {Histogram maximum likelihood estimator of a
periodic function in the multiplicative intensity model}. {\it
Statistics and Decisions} 6: 79-88.
\bibitem{leskow}
Leśkow J., Różański R. (1989) {Histogram maximum likelihood
estimator in the multiplicative intensity model}. {\it Stochastic
Processes and their Applications} 31: 151-159.
\bibitem{wronka}
Leśkow J., Wronka C. (2004) {Bootstrap resampling in the analysis
of time series}. {\it Springer-Verlang, Heidelberg-Berlin}, pp.
267-274.
\bibitem{loader}
Loader C. (1993) {Nonparemetric regression, confidence bands and
bias correction}. {\it Computing Science and Statistics}:
Proceedings of the 25th Symposium on the Interface, 131-136.
\bibitem{sneth}
Snethlage M. (1999) { Is bootstrap really helpful in point process
statistics?} {\it Metrika} 49: 245-255.


\end{thebibliography}
\end{document}